\documentclass[12pt,a4paper]{article}
\usepackage[cp1251]{inputenc} 
\usepackage[russian]{babel}
\usepackage{amsfonts}

\newcommand{\di}{\displaystyle}

\newcommand{\al}{\alpha}

\newcommand{\ga}{\gamma}
\newcommand{\de}{\delta}
\newcommand{\la}{\lambda}
\newcommand{\om}{\omega}
\newcommand{\ee}{\varepsilon}

\newcommand{\iy}{\infty}

\begin{document}

\begin{center}
{\large\bf
On Inverse Spectral Problems for Second Order Integro-differential
Operators}\\[0.1cm]
{\bf V.A.\,Yurko} \\[0.1cm]
\end{center}

\thispagestyle{empty}

\noindent {\bf Abstract.} Inverse spectral problems are studied for the second order integro-differential operators on a finite interval.
Properties of spectral characteristic are established, and the uniqueness theorem is proved for this class of inverse problems.



{\bf 1. } Inverse spectral problems consist in recovering operators
from their spectral characteristics. Such problems often appear in
mathematics, mechanics, physics, electronics, geophysics and other
branches of natural sciences and engineering. The greatest success
in the inverse problem theory has been achieved for the Sturm-Liouville
operator (see, e.g., [1-3]) and afterwards for higher-order differential
operators [4-6] and other classes of differential operators.

For integro-differential and other classes of nonlocal operators inverse
problems are more difficult for investigation, and the main classical methods
(transformation operator method and the method of spectral mappings [1-6])
either are not applicable to them or require essential modifications, and
for such operators the general inverse problem theory does not exist.
At the same time, nonlocal and, in particular, integro-differential operators
are of great interest, because they have many applications (se, e.g., [7]).
We note that some aspects of inverse problems for integro-differential
operators were studied in [8-10] and other works. In the present paper
we study inverse spectral problem for one class of second order
integro-differential operators on a finite interval. Properties of
spectral characteristic are established, and the uniqueness theorem
is proved for this class of inverse problems.

\bigskip
{\bf 2. } Consider the integro-differential equation
$$
\ell y:=-y''(x)+q(x)y(x)+\int_0^x M(x,t)y(t)\,dt=\la y(x),\;x\in[0,\pi],    \eqno(1)
$$
where $q(x)$ and $M(x,t)$ are integrable complex-valued functions. Let
$C(x,\la)$ and $S(x,\la)$ be solutions of Eq. (1) with the initial conditions
$$
C(0,\la)=S'(0,\la)=1,\quad C'(0,\la)=S(0,\la)=0.
$$
For each fixed $x\in[0,\pi],$ the functions $C^{(\nu)}(x,\la)$ and
$S^{(\nu)}(x,\la),$ $\nu=0,1,$ are entire in $\la$ of order $1/2.$
Denote $\Delta_1(\la):=S(\pi,\la),$ $\Delta_2(\la):=C(\pi,\la).$
Zeros $\{\la_{nk}\}_{n\ge 1}$ of the entire function $\Delta_k(\la)$
coincide with the eigenvalues of the boundary value problem $L_k=L_k(M,q)$
for Eq. (1) with the conditions $y^{(k-1)}(0)=y(\pi)=0.$ The function
$\Delta_k(\la)$ is called the characteristic function for $L_k$.

Let $\Phi(x,\la)$ be the solution of Eq. (1) under the conditions
$\Phi(0,\la)=1,$ $\Phi(\pi,\la)=0.$ Denote $N(\la):=\Phi'(0,\la).$ Then
$$
\Phi(x,\la)=C(x,\la)+N(\la)S(x,\la),
\quad N(\la)=-\Delta_2(\la)/\Delta_1(\la).                                 \eqno(2)
$$
The function $N(\la)$ is called the Weyl-type function for $\ell.$ It
follows from (2) that the function $N(\la)$ is meromorphic in $\la$
with poles $\{\la_{n1}\}_{n\ge 1}$ and zeros $\{\la_{n2}\}_{n\ge 1}$.
Let $M(x,t)$ be known a priori.
The inverse problem is formulated as follows.

\smallskip
{\bf Inverse problem 1. }{\it Given $N(\la),$ construct $q(x).$}

\smallskip
This inverse problem is an analogue of the classical inverse problem
of recovering the Sturm-Liouville operator from the given Weyl function [3].

\medskip
{\bf 3. } Let $\la=\rho^2.$ By the well-known method (see, e.g., [3])
we have for $|\rho|\to\iy$:
$$
S(x,\la)=\frac{\sin\rho x}{\rho}-\frac{\cos\rho x}{2\rho^2}\int_0^x q(t)\,dt
+o\Big(\frac{1}{\rho}\Big),\;
C(x,\la)=\cos\rho x+\frac{\sin\rho x}{2\rho}\int_0^x q(t)\,dt
+o\Big(\frac{1}{\rho}\Big),
$$
and consequently,
$$
\Delta_1(\la)=\frac{\sin\rho\pi}{\rho}-
\frac{\om\cos\rho\pi}{2\rho^2}+o\Big(\frac{1}{\rho}\Big),\;
\Delta_2(\la)=
\cos\rho\pi+\frac{\om\sin\rho\pi}{2\rho}+o\Big(\frac{1}{\rho}\Big),      \eqno(3)
$$
where $\om:=\di\int_0^\pi q(t)\,dt.$
Using (3), by standard calculations [3] we obtain
$$
\rho_{n1}:=\sqrt{\la_{n1}}=n+\frac{\om}{2n}+o\Big(\frac{1}{n}\Big),\;
\rho_{n2}:=\sqrt{\la_{n2}}=n-\frac{1}{2}+\frac{\om}{2n}+o\Big(\frac{1}{n}\Big).
$$
Moreover, the specification of $\{\la_{nk}\}_{n\ge 1}$ uniquely determines
the characteristic function $\Delta_k(\la)$ by the formulae [3]:
$$
\Delta_1(\la)=\pi\prod_{n=1}^{\iy}\frac{\la_{n1}-\la}{n^2}\,,\;
\Delta_2(\la)=\prod_{n=1}^{\iy}\frac{\la_{n2}-\la}{(n-1/2)^2}\,.        \eqno(4)
$$
Taking (2) and (4) into account we conclude that Inverse problem 1
is equivalent to the following Borg-type inverse problem.

\smallskip
{\bf Inverse problem 2. }{\it Given two spectra
$\{\la_{nk}\}_{n\ge 1},\; k=1,2,$ construct $q(x).$}

\smallskip
Denote
$$
\ell^{*}z:=-z''(x)+q(x)z(x)+\int_x^{\pi} M(t,x)z(t)\,dt.
$$
Obviously,
$$
\int_0^\pi \ell y(x)\cdot z(x)=
\Big|_0^\pi (yz'-y'z)+\int_0^\pi y(x)\cdot\ell^{*}z(x).                \eqno(5)
$$
If $y(x,\la)$ and $z(x,\mu)$ are solutions of the equations
$\ell y=\la y$ and $\ell^{*} z=\mu z,$ respectively, then (5) yields
$$
(\la-\mu)\int_0^\pi y(x,\la)z(x,\mu)\,dx=\Big|_0^\pi (yz'-y'z).       \eqno(6)
$$
Let $C_{*}(x,\la)$ and $S_{*}(x,\la)$ be solutions of the equation
$$
\ell^{*} z=\la z                                                      \eqno(7)
$$
with the initial conditions
$$
C_{*}(\pi,\la)=-S_{*}'(\pi,\la)=1,\quad C_{*}'(\pi,\la)=S_{*}(\pi,\la)=0.
$$
Denote $\Delta^{*}_1(\la):=S_{*}(0,\la),$ $\Delta^{*}_2(\la):=
-S_{*}'(0,\la).$ It follows from (6) with $\mu=\la$ that
$$
S(\pi,\la)\equiv S_{*}(0,\la),\; S'(\pi,\la)\equiv C_{*}(0,\la),\;
C(\pi,\la)\equiv -S'_{*}(0,\la),\; C'(\pi,\la)\equiv -C_{*}(0,\la),  \eqno(8)
$$
hence,
$$
\Delta^{*}_1(\la)\equiv\Delta_1(\la),\quad
\Delta^{*}_2(\la)\equiv\Delta_2(\la).
$$
Let $\Phi_{*}(x,\la)$ be the solution of Eq. (7) under the conditions
$\Phi_{*}(0,\la)=1,$ $\Phi_{*}(\pi,\la)=0.$ Denote $N^{*}(\la):=
\Phi'_{*}(0,\la).$ Then
$$
\Phi_{*}(x,\la)=\frac{S_{*}(x,\la)}{S_{*}(0,\la)},\quad
N^{*}(\la)=\frac{S'_{*}(0,\la)}{S_{*}(0,\la)}\,.
$$
Together with (8) this yields $N^{*}(\la)\equiv N(\la).$

\smallskip
It is known (see, e.g., [11]) that there exists a fundamental
system of solutions $\{y_1(x,\rho), y_2(x,\rho)\},$ $\mbox{Im}\,\rho
\ge 0,\,x\in[0,\pi]$ for Eq. (1) such that for $|\rho|\to\iy,\;\nu=0,1$:
$$
y_1^{(\nu)}(x,\rho)=(i\rho)^{\nu}\exp(i\rho x)+O(\rho^{-1}),\;
y_2^{(\nu)}(x,\rho)=(-i\rho)^{\nu}\exp(-i\rho x)(1+O(\rho^{-1})).
$$
Similarly, there exists a fundamental system of solutions
$\{z_1(x,\rho), z_2(x,\rho)\},\,\mbox{Im}\,\rho\ge 0,\,x\in[0,\pi]$
for Eq. (7) such that for $|\rho|\to\iy,\;\nu=0,1$:
$$
z_1^{(\nu)}(x,\rho)=(i\rho)^{\nu}\exp(i\rho x)(1+O(\rho^{-1})),\;
z_2^{(\nu)}(x,\rho)=(-i\rho)^{\nu}\exp(-i\rho x)+O(\rho^{-1}\exp(-i\rho\pi)).
$$
Fix $\de,\ee\in(0,\pi/2).$ Denote $Q:=\{\rho:\;\arg\rho\in[\de,\pi-\de]\}.$
Using these fundamental systems of solutions we obtain the following
asymptotics for $\rho\in Q,\,|\rho|\to\iy,$ uniformly in $x\in[0,\pi-\ee]$:
$$
\Phi(x,\la)=\exp(i\rho x)(1+O(\rho^{-1}))+O(\rho^{-1}),
\quad \Phi_{*}(x,\la)=\exp(i\rho x)(1+O(\rho^{-1})).                   \eqno(9)
$$

\smallskip
{\bf 4. } In this section we provide an algorithm for the solution of
Inverse problem 1. For this purpose together with $L_k$ we consider
the boundary value problems $\tilde L_k:=L_k(M,\tilde q)$ of the same
form but with a different potential $\tilde q(x).$ We agree that
everywhere below if a certain symbol $\al$ denotes an object related
to $L_k$, then $\tilde\al$ will denote the analogous object related
to $\tilde L_k$, and $\hat\al:=\al-\tilde\al.$

\smallskip
{\bf Lemma 1. }{\it Let $M(x,t)\equiv\tilde M(x,t).$ Then}
$$
\int_0^\pi\hat q(x)\Phi(x,\la)
\tilde\Phi_{*}(x,\la)\,dx\equiv -\hat N(\la),                       \eqno(10)
$$
where $\hat q(x)=q(x)-\tilde q(x),\;\hat N(\la)=N(\la)-\tilde N(\la).$

\smallskip
{\it Proof. } One has
$$
-\Phi''(x,\la)+q(x)\Phi(x,\la)+\int_0^x M(x,t)\Phi(t,\la)\,dt
=\la\Phi(x,\la),
$$
$$
-\tilde\Phi''_{*}(x,\la)+\tilde q(x)\tilde\Phi_{*}(x,\la)+
\int_x^\pi M(t,x)\tilde\Phi_{*}(t,\la)\,dt=\la\tilde\Phi_{*}(x,\la).
$$
We multiply the first relation by $\tilde\Phi_{*}(x,\la),$
then subtract the second relation multiplying by
$\Phi(x,\la)$ and integrate with respect to $x:$
$$
\int_0^\pi\hat q(x)\Phi(x,\la)\tilde\Phi_{*}(x,\la)\,dx+
\int_0^\pi \tilde\Phi_{*}(x,\la)\,dx\int_0^x M(x,t)\Phi(t,\la)\,dt
$$
$$
-\int_0^\pi \Phi(x,\la)\,dx\int_x^\pi M(t,x)\tilde\Phi_{*}(t,\la)\,dt
=\Big|_0^\pi \Big(\Phi'(x,\la)\tilde\Phi_{*}(x,\la)
-\Phi(x,\la)\tilde\Phi'_{*}(x,\la)\Big).
$$
Taking the relations $\Phi(0,\la)=1,\,\Phi(\pi,\la)=0,
\,\tilde\Phi_{*}(0,\la)=1,\,\tilde\Phi_{*}(\pi,\la)=0,\,
N(\la)=\Phi'(0,\la),$ $\tilde N^{*}(\la)=\tilde\Phi'_{*}(0,\la)$
and $\tilde N^{*}(\la)=\tilde N(\la)$ into account
we arrive at (10).

\smallskip
{\bf Lemma 2. }{\it Let}
$$
r(x)=\di\frac{x^k}{k!}\Big(\gamma+p(x)\Big), \quad
H(x,\rho)=\exp(2i\rho x)\Big(1+\di\frac{\xi(x,\rho)}{\rho}\Big)
+\exp(i\rho x)\frac{\eta(x,\rho)}{\rho},\; x\in[0,\pi],
$$
{\it where $p(x)\in C[0,\pi],\;p(0)=0,$ and where the functions
$\xi(x,\rho), \eta(x,\rho)$ are continuous and bounded for
$x\in[0,\pi],\;\rho\in Q,\;|\rho|\ge\rho^*.$ Then for}
$|\rho|\to\infty,\;\rho\in Q,$
$$
\di\int_0^\pi r(x)H(x,\rho)\,dx=\di\frac{1}{(-2i\rho)^{k+1}}(\ga+o(1)).
$$

\smallskip
{\it Proof.} We calculate
$$
(-2i\rho)^{k+1}\di\int_0^\pi r(x)H(x,\rho)\,dx
=I_1(\rho)+I_2(\rho)+I_3(\rho)+I_4(\rho),
$$
where
$$
I_1(\rho)=\gamma(-2i\rho)^{k+1}\di\int_0^\pi
\di\frac{x^k}{k!}\exp(2i\rho x)\,dx,
$$
$$
I_2(\rho)=(-2i\rho)^{k+1}\di\int_0^\pi \di\frac{x^k}{k!}p(x)
\exp(2i\rho x)\,dx,
$$
$$
I_3(\rho)=(-2i)^{k+1}\rho^k\di\int_0^\pi r(x)\exp(2i\rho x)
\xi(x,\rho)\,dx,
$$
$$
I_4(\rho)=(-2i)^{k+1}\rho^k\di\int_0^\pi r(x)\exp(i\rho x)
\eta(x,\rho)\,dx,
$$
Since
$$
\di\int_0^{\infty} \di\frac{x^k}{k!}\exp(2i\rho x)\,dx=
\di\frac{1}{(-2i\rho)^{k+1}}, \quad \rho\in Q,
$$
it follows that $I_1(\rho)-\gamma\to 0$ as $|\rho|\to\iy,\;\rho\in Q.$
If $\rho\in Q,$ then there exists $\varepsilon_0 >0$ such that
$$
|Im\,\rho|\ge\varepsilon_0|\rho|\;\mbox{ for }\;\rho\in Q.     \eqno(11)
$$
Take $\varepsilon > 0$ and choose $\delta=\delta(\varepsilon)$
such that for $x\in[0,\delta],\quad |p(x)|<\di\frac{\varepsilon}{2}
\varepsilon_0^{k+1},$ where $\varepsilon_0$ is defined in (11).
Then, using (11) we infer
$$
|I_2(\rho)|<\di\frac{\varepsilon}{2}(2|\rho|\varepsilon_0)^{k+1}
\di\int_0^{\delta} \di\frac{x^k}{k!}\exp(-2\varepsilon_0|\rho| x)\,dx+
(2|\rho|)^{k+1}\di\int_{\delta}^{\pi} \di\frac{x^k}{k!}|p(x)|
\exp(-2\varepsilon_0|\rho| x)\,dx
$$
$$
<\di\frac{\varepsilon}{2}+(2|\rho|)^{k+1}\exp(-2\varepsilon_0|\rho|\delta)
\di\int_{0}^{\pi-\delta} \di\frac{(x+\delta)^k}{k!}|p(x+\delta)|
\exp(-2\varepsilon_0|\rho| x)\,dx.
$$
By arbitrariness of $\varepsilon$ we obtain that $I_2(\rho)\to 0$ for
$|\rho|\to\infty,\;\rho\in Q.$

Since $|(\gamma+p(x))\xi(x,\rho)|<C,$ then for $|\rho|\to\infty,\;
\rho\in Q,$
$$
|I_3(\rho)|<C|\rho|^k\di\int_0^{\pi} \di\frac{x^k}{k!}
\exp(-2\varepsilon_0|\rho| x)\,dx \le
\di\frac{C}{|\rho|\varepsilon_0^{k+1}},
$$
hence $I_3(\rho)\to 0$ for $|\rho|\to\infty,\;\rho\in Q.$
Similarly, one gets $I_4(\rho)\to 0$ for $|\rho|\to\infty,\;\rho\in Q.$
$\hfill\Box$

\smallskip
For simplicity, we will assume that $q(x)$ is analytic on $[0,\pi].$
Suppose that for a certain fixed $k\ge 0$ the Taylor coefficients
$q_j:=q^{(j)}(0), \; j=\overline{0,k-1},$ have been already found.
Let us choose a model potential $\tilde q(x)$ such that
the first $k$ Taylor coefficients of $q$ and $\tilde q$ coincide,
i.e. $\tilde q_j=q_j, \; j=\overline{0,k-1}.$ Then, using (9)-(10)
and Lemma 2, we can calculate the next Taylor coefficient
$q_k=q^{(k)}(0).$ Namely, the following assertion is valid.

\smallskip
{\bf Lemma 3. }{\it Fix $k.$ Let the functions $q(x)$ and
$\tilde q(x)$ be analytic for $x \in [0,\pi],$ with $\hat q_j:=
q_j-\tilde q_j=0$ for $j=\overline{0,k-1}.$ Then}
$$
\hat q_k=-\di\lim_{|\rho|\to\infty
\atop\rho\in Q}\;(-2i\rho)^{k+1}\hat N(\lambda).               \eqno(12)
$$

\smallskip
Thus, we arrive at the following algorithm for the solution of
Inverse Problem 1.

\smallskip
{\bf Algorithm 1. }{\it Let the Weyl-type function $N(\la)$
be given. Then:

(i) We calculate $q_k=q^{(k)}(0),\; k\ge 0.$ For this purpose we
successively perform the following operations for $k=0,1,2,\ldots:$
We construct a model potential $\tilde q(x)$ such that
$\tilde q_j=q_j,\; j=\overline{0,k-1}$ and arbitrary in the rest,
and we calculate $q_k=q^{(k)}(0)$ by (12).

(ii) We construct the function $q(x)$ by the formula
$$
q(x)=\di\sum_{k=0}^{\infty}q_k\di\frac{x^k}{k!}, \quad 0<x<R,
$$
where
$$
R=\Big(\overline{\di\lim_{k\to\infty}}\,
\Big(\di\frac{|q_k|}{k!}\Big)^{1/k}\Big)^{-1}.
$$
If $R<\pi$ then for $R<x<\pi$ the function $q(x)$ is constructed
by analytic continuation.}

\bigskip
{\bf Acknowledgment.} This work was supported in part by Grant 1.1660.2017/PCh
of the Russian Ministry of Education and Science and by Grants 16-01-00015,
17-51-53180 of Russian Foundation for Basic Research.

\begin{center}
{\bf REFERENCES}
\end{center}
\begin{enumerate}
\item[{[1]}] Marchenko V.A., Sturm-Liouville operators and their applications.
     "Naukova Dumka",  Kiev, 1977;  English  transl., Birkh\"auser, 1986.
\item[{[2]}] Levitan B.M., Inverse Sturm-Liouville problems. Nauka,
     Moscow, 1984; English transl., VNU Sci.Press, Utrecht, 1987.
\item[{[3]}] Freiling G. and Yurko V.A., Inverse Sturm-Liouville Problems
     and their Applications. NOVA Science Publishers, New York, 2001.
\item[{[4]}] Beals R., Deift P. and Tomei C.,  Direct and Inverse Scattering
     on the Line, Math. Surveys and Monographs, v.28. Amer. Math. Soc.
     Providence: RI, 1988.
\item[{[5]}] Yurko V.A. Method of Spectral Mappings in the Inverse Problem
     Theory. Inverse and Ill-posed Problems Series. VSP, Utrecht, 2002.
\item[{[6]}] Yurko V.A. Inverse Spectral Problems for Differential
     Operators and their Applications. Gordon and Breach, Amsterdam, 2000.
\item[{[7]}] Lakshmikantham V. and Rama Mohana Rao M. Theory of
     integro-differential equations. Stability and Control: Theory and
     Applications, vol.1, Gordon and Breach, Singapore, 1995.
\item[{[8]}] Yurko V.A., An inverse problem for integro-differential
     operators. Matem. zametki, 50, no.5 (1991), 134-146 (Russian);
     English transl. in Math. Notes, 50, no.5-6 (1991), 1188-1197.
\item[{[9]}] Kuryshova Yu. An inverse spectral problem for differential
     operators with integral delay. Tamkang J. Math. 42, no.3 (2011), 295-303.
\item[{[10]}] Buterin S.A. On the reconstruction of the convolution
     perturbation of the Sturm-Liouville operator from the spectrum,
     Differential Equations 46, no.1 (2010), 150--154.
\item[{[11]}] Hromov A.P. On generating functions of Volterra operators.
     Math. USSR Sbornik 31, no.3 (1997), 409-432.
\end{enumerate}

\end{document}